\documentclass{article}
\usepackage{amsmath,amscd,amssymb,amsfonts,amsthm}   
\usepackage[all]{xy}   %for commutative diagrams
\UseComputerModernTips %to improve the arrows in commutative diagrams
\CompileMatrices
\theoremstyle{plain}
	\newtheorem{theorem}{Theorem}[section]
        \newtheorem{maintheorem}{Theorem}

	\newtheorem{notation}[theorem]{Notation}
	\newtheorem{example}[theorem]{Example}
	\newtheorem{lemma}[theorem]{Lemma}
	\newtheorem{proposition}[theorem]{Proposition}
	\newtheorem{corollary}[theorem]{Corollary}

\theoremstyle{definition}
	\newtheorem{definition}[theorem]{Definition}

\DeclareMathOperator{\Ker}{Ker}
\DeclareMathOperator{\Image}{Image}

\DeclareMathOperator{\Aut}{Aut}
\DeclareMathOperator{\End}{End}
\DeclareMathOperator{\GL}{GL}
\DeclareMathOperator{\id}{id}
\def\ab{\mathrm{ab}}
\def\Z{\mathbb Z}
\def\N{\mathbb N}
\def\Q{\mathbb Q}
\def\functor{\underline{\phantom{A}}} 
\def\talpha{\widetilde{\alpha}}
\def\tbeta{\widetilde{\beta}}
\def\Nov{\mathrm{Nov}}
\def\mfa{{\mathfrak{a}}}
\def\mfb{{\mathfrak{b}}}
\def\tEnd{\widetilde{\End}}
\def\tK{\widetilde{K}}
\DeclareMathOperator{\rad}{rad}
\DeclareMathOperator{\tAuto}{\underline{\widetilde{\Aut}}}
\def\cy#1{{\mathcal #1}}
\begin{document}
\begin{center}
{\Large Whitehead Groups of Localizations and the Endomorphism Class
Group} \\[2ex]
{\large Desmond Sheiham}{\def\thefootnote{}\footnote{Submitted September
24th 2002. In revised form June 25th 2003.}\footnote{2000 Mathematics
    Subject Classification 18F25; 16S10; 16S34; 37C27.}\addtocounter{footnote}{-1}}
\end{center}
%\begin{frontmatter}
%\title{Whitehead Groups of Localizations and the Endomorphism Class
%Group}
%\author{Desmond Sheiham}
%\address{Department of Mathematics \\ University of California,
%  Riverside \\ CA 92521, USA \\ des@sheiham.com}
\begin{abstract}
We compute the Whitehead groups of the associative rings in a class
which includes (twisted) formal power series rings and the augmentation
localizations of group rings and polynomial rings.

For any associative ring $A$, we obtain an invariant of a pair
$(P,\alpha)$ where $P$ is a finitely generated projective $A$-module 
and $\alpha:P\to P$ is an endomorphism. This
invariant determines $(P,\alpha)$ up to extensions, 
yielding a computation of the (reduced) endomorphism class group
$\tEnd_0(A)$.

We also refine the analysis by Pajitnov and Ranicki of the
Whitehead group of the Novikov ring, a computation which Pajitnov has
used in work on circle-valued Morse theory.
\end{abstract}
%\end{frontmatter}
\section{Introduction}
\subsection{Endomorphisms}
The characteristic polynomial of an endomorphism of a vector space
determines the endomorphism uniquely `up to choices of extension'.
To make such a statement precise, one 
makes the following definition (Almkvist\cite{Alm73,Alm74,Alm78},
Grayson~\cite{Gra77}) which we discussed in~\cite{She01}:
\begin{definition}
\label{end_class_gp_free}
Let $A$ be an associative ring. The reduced endomorphism class group
$\tEnd_0(A)$ is the abelian group with one generator for each
isomorphism class of pairs $[A^n,\alpha]$ and relations: 
\begin{itemize}
\item
$[A^n,\alpha]+[A^{n''},\alpha''] = [A^{n'},\alpha']$ if there is an
exact sequence 
\begin{equation*}
0\to (A^n,\alpha)\to (A^{n'},\alpha') \to
(A^{n''},\alpha'')\to 0.
\end{equation*}
\item $[A,0]=0$.
\end{itemize}
\end{definition}
Almkvist proved~\cite{Alm74} that if $A$ is commutative then the
characteristic polynomial induces an isomorphism 
\begin{align}
\label{Almkvist_result}
\tEnd_0(A) &\cong \left\{ \frac{1+a_1x
+\cdots+a_nx^n}{1+b_1x+\cdots + b_mx^m} \bigm{|} a_i,b_i\in A\right\} \\
[A^n,\alpha] & \mapsto \det(1-\alpha x) \notag
\end{align}
A goal of the present paper is to obtain an analogous
statement for arbitrary associative rings. 

We first reformulate the right-hand side
of~(\ref{Almkvist_result}). Let $\epsilon: A[x]\to A;x\mapsto 0$ denote
augmentation and let $P\subset A[x]$ denote the set of polynomials $p$
such that $\epsilon(p)$ is invertible. There is a canonical factorization
of $\epsilon$ through the ring of fractions $P^{-1}A[x]$ 
\begin{equation*}
A[x]\rightarrowtail P^{-1}A[x] \xrightarrow{\epsilon_P} A
\end{equation*}
and Almkvist's result says that for any commutative ring $A$ there is
an isomorphism 
$\tEnd_0(A)\cong\epsilon_P^{-1}(1)$.

Dropping the assumption that $A$ is commutative, let $A[x]$ denote the
polynomial ring in a central indeterminate $x$.
It is advantageous to make invertible not only elements but
matrices in $A[x]$. Let $\Sigma$ denote the set of
matrices $\sigma$ with entries in $A[x]$ such that $\epsilon(\sigma)$
is invertible. There exists a (formal) localization $A[x]\to\Sigma^{-1}A[x]$
(Cohn~\cite{Coh71}, Schofield~\cite[Ch4]{Scho85})  
which is not in general a ring of fractions but has the following
properties:
\begin{enumerate}
\item\label{invertiblizes} Every matrix $\sigma\in\Sigma$ is invertible over
$\Sigma^{-1}A[x]$.
\item The map $A[x]\to \Sigma^{-1}A[x]$ is universal with respect
to property~\ref{invertiblizes}. In other words, any map $A[x]\to B$ which
makes all the matrices in $\Sigma$ invertible factors uniquely through
$\Sigma^{-1}A[x]$. 
\end{enumerate}
In particular $\epsilon:A[x]\to A$ is the composite
$\displaystyle{A[x]\rightarrowtail \Sigma^{-1}A[x] \xrightarrow{\epsilon_\Sigma} A}$.
\begin{maintheorem}
\label{Endomorphism_class_group}
Let $A$ be an associative ring. There is an isomorphism
\begin{align*}
\tEnd_0(A) &\cong \epsilon_\Sigma^{-1}(1)/C \\
[A^n,\alpha] &\mapsto D(1-\alpha x)
\end{align*}
where $C$ is the subgroup generated by commutators:
\begin{equation*}
\{(1+ab)(1+ba)^{-1} \mid a,b\in A[x],~\epsilon(ab)=\epsilon(ba)=0\}.
\end{equation*}
\end{maintheorem}
The symbol $D$ is defined in section~\ref{section:proof} below
(definition~\ref{Define_D}) and is
analogous to the Dieudonn\'e determinant. 
In the commutative case a matrix in $\Sigma$ is invertible if and only
if its determinant is invertible, so
$\Sigma^{-1}A[x]=P^{-1}A[x]$. Moreover, $C$ is trivial and $D$ is the
traditional determinant so theorem~\ref{Endomorphism_class_group} is a
generalization of Almkvist's identity~(\ref{Almkvist_result}) above.

Suppose $A$ is non-commutative. Now
\begin{equation*}
[\epsilon_\Sigma^{-1}(1),\epsilon_\Sigma^{-1}(1)]~\subset~\left[(\Sigma^{-1}A[x])^\bullet,(\Sigma^{-1}A[x])^\bullet\right]\cap\epsilon_\Sigma^{-1}(1)~\subset~C
\end{equation*}
where $(\Sigma^{-1}A[x])^\bullet$ denotes the group of units in
$\Sigma^{-1}A[x]$ but neither of these inclusions is an equality in
general and
$[\epsilon_\Sigma^{-1}(1),\epsilon_\Sigma^{-1}(1)]\subset C$ is never
an equality (see section~\ref{section:study_C}).
%and we emphasize that it is usually much larger than the preimage of
%$1\in A$ %in the {\em ring} $\Sigma^{-1}A[x]/I$ where $I$ is the
%two-sided ideal 
%generated by commutators $\{ab-ba\mid \epsilon(ab)=\epsilon(ba)=0\}$.

Our proof of theorem~\ref{Endomorphism_class_group} uses a 
result of A.Ranicki which we state next. Recall that for an
arbitrary ring $A$, the group $K_1(A)=\GL(A)^{\ab}$ is the 
abelianization of the group of invertible square matrices of arbitrary size.
Ranicki established~\cite[Prop10.21]{Ran98} an isomorphism
\begin{equation*}
K_1(A)\oplus \tEnd_0(A) \xrightarrow{\cong} K_1(\Sigma^{-1}A[x])
\end{equation*}
which is the canonical inclusion of $K_1(A)$ and is defined on 
$\tEnd_0(A)$ by
\begin{equation*}
[A^n,\alpha] \mapsto [1-\alpha x]\in
K_1(\Sigma^{-1}A[x]).
\end{equation*} We prove here that there is an isomorphism $D$
between the image of $\tEnd_0(A)$ in $K_1(\Sigma^{-1}A[x])$ and
$\epsilon_\Sigma^{-1}(1)/C$.
\subsection{Local Augmentations}
Our main result is more general and concerns a class of ring homomorphisms
$\epsilon:B\to A$ which we call `local augmentations'. The word 
`augmentation' just means split surjection or in other words
retraction in the category of rings. By local we 
mean that a square matrix $\alpha$ with entries in $B$ is invertible
if $\epsilon(\alpha)$ is invertible. 

Any augmentation $\epsilon:B\to A$ can be made a local augmentation
$\Sigma^{-1}B\to A$ by adjoining a formal inverse to
every square matrix $\sigma$ with entries in $B$ such that
$\epsilon(\sigma)$ is invertible (see \cite[Lemma 3.1]{She01} and
lemmas~\ref{localize_morphism} and~\ref{loc_and_complete} below).
In particular, the map $\epsilon_\Sigma:\Sigma^{-1}A[x]\to A$ above is a
local augmentation. To give another example, if $A[[x]]$ denotes
the ring of formal power series in a central indeterminate $x$ with
coefficients in $A$ then $A[[x]]\to A;x\to 0$ is a local
augmentation. 
\begin{maintheorem}[Main Theorem]
\label{identify_K1}
If $\epsilon:B\to A$ is a local augmentation then there is a canonical
isomorphism
\begin{equation*}
K_1(B)\cong K_1(A) \oplus \frac{\epsilon^{-1}(1)}{C}.
\end{equation*}
where $C$ is the subgroup generated by commutators:
\begin{equation*}
\{(1+ab)(1+ba)^{-1} \mid a,b\in B,~\epsilon(ab)=\epsilon(ba)=0\}
\end{equation*}
\end{maintheorem}
\noindent In particular, we may apply theorem~\ref{identify_K1} to
augmentation localizations of group rings. Suppose $\pi$ is a group and
$A\pi=A\otimes_\Z\Z\pi$ is the corresponding group ring with coefficients in an
associative ring $A$. Let $\epsilon:A\pi \to A$ be the
augmentation, defined by $\epsilon(g)=1$ for $g\in\pi$ and
$\epsilon(a)=a$ for $a\in A$. Let 
$\Sigma$ denote the set of matrices $\sigma$ with entries in $A\pi$ such that
$\epsilon(\sigma)$ is invertible; $\epsilon$ can be written as the
composite $A\pi\to \Sigma^{-1}A\pi \xrightarrow{\epsilon_\Sigma} A$.
\begin{corollary}
\label{group_ring_corollary}
\qquad$\displaystyle{K_1(\Sigma^{-1}A\pi) \cong K_1(A) \oplus
\frac{\epsilon_\Sigma^{-1}(1)}{C}}$. \hspace*{\fill}\qed
\end{corollary}
Corollary~\ref{group_ring_corollary} will be applied in a subsequent
paper, with $A=\Z$, to study Reidemeister torsion 
of {\em homology} equivalences between finite CW-complexes.

We may also apply theorem~\ref{identify_K1} to the ring of formal power
series in a central indeterminate:
\begin{equation}
\label{untwisted_power_series}
K_1(A[[x]])\cong K_1(A)\oplus \frac{1+A[[x]]x}{C}
\end{equation}
It may be useful, if one is studying endomorphisms via
theorem~\ref{Endomorphism_class_group}, to pass from the localization
$\Sigma^{-1}A[x]$ to the completion $A[[x]]$. The universal property
of localization provides a canonical map $\gamma:\Sigma^{-1}A[x]\to
A[[x]]$ - see lemma~\ref{loc_and_complete}d) - but neither 
$\gamma$ nor the induced map
\begin{equation*}
\frac{\epsilon_\Sigma^{-1}(1)}{C}\to
\frac{1+A[[x]]x}{C}
\end{equation*}
 is an injection in general~\cite{She01}.

Generalizing~(\ref{untwisted_power_series}), suppose $A$ is an
associative ring and $\xi:X\to\Aut(A)$ assigns a ring automorphism to
each element of a set $X$ of indeterminates. 
Let $A_\xi\langle\langle X\rangle\rangle$ denote the (twisted) power
series ring whose elements are infinite formal sums, with one term for each
word in the alphabet $X$. One may also impose relations, such as
commutativity of the indeterminates, if compatible with $\xi$; see
example~\ref{twisted_power_series_example} below.
\begin{corollary}
\label{power_series_corollary}
\qquad$\displaystyle{K_1(A_\xi\langle\langle X\rangle\rangle)\cong K_1(A)\oplus
\frac{1+A_\xi\langle\langle X\rangle\rangle X}{C}}$. \hspace*{\fill}\qed
\end{corollary}
We remark that commutators of the form $(1+ab)(1+ba)^{-1}$ have
appeared in earlier computations related to Whitehead groups
(e.g. \cite[p269]{Bas68},~\cite{MenMon84},~\cite{Sil81}). 
Whitehead groups of universal localizations appear in the $K$-theory
exact sequences of Schofield~\cite{Scho85} and Neeman and
Ranicki~\cite{NeeRan01} and certain Whitehead groups
of localizations are computed in papers of Revesz~\cite{Rev83}, Ara,
Goodearl and Pardo~\cite{AGP02} and Ara~\cite{Ara03}.

\subsection{The Novikov Ring}
Corollary~\ref{power_series_corollary} refines a
computation by Pajitnov and Ranicki~\cite{PajRan00} which was
motivated by work of Pajitnov~\cite{Paj00} on circle-valued Morse theory and
gradient flow on manifolds. We briefly outline this application;
background references include~\cite{Bot88,Mil63,Mil65,Nov81,Nov82,Paj95,Paj99,Ran01}.

Let $A_\xi((z))=A_\xi[[z]][z^{-1}]$ denote the (twisted) Novikov ring
whose elements are power series with finitely many negative powers of
$z$ but, in general, infinitely many positive powers of $z$.
Let $M$ be a (closed compact) manifold, $f:M\to S^1$ a Morse
map and $v:M\to \tau_M$ a generic choice of vector field which is
`gradient-like' for $f$. One would like to describe the closed orbits
of the associated flow on $M$.

Now if $f$ induces a surjection $\pi_1(M)\to \pi_1(S^1)=\Z$ then
\begin{equation*}
\pi_1(M)\cong \pi\rtimes_\xi\Z = \pi \rtimes_\xi \{z^n\mid n\in\Z\}
\end{equation*}
for some homomorphism $\xi:\Z\to\Aut(\pi)$ and the Novikov ring $\Z\pi_\xi((z))$ is a completion of the
group ring $\Z\pi_1(M)$. Information about
the closed orbits of the gradient flow $v$ is encoded in the torsion
\begin{equation*}
\tau(\phi_v)\in \frac{K_1(\Z\pi_\xi((z)))}{\Image(\pm\pi_1M)}
\end{equation*}
of a canonical chain equivalence (Novikov~\cite{Nov81,Nov82}, Pajitnov~\cite{Paj95,Paj99,Paj00}) 
\begin{equation*}
\phi_v:\Z\pi_\xi((z))\otimes_{\Z{\pi_1(M)}}C^\Delta \to C^\Nov(v).
\end{equation*}
The symbol $C^\Delta$ denotes the chain complex over
$\Z\pi_1(M)$ associated to a smooth triangulation of $M$. The `Novikov
complex' $C^\Nov(v)$ is a 
finitely generated free chain complex over $\Z\pi_\xi((z))$
with one basis element in $C^\Nov_i$ for each critical point of $f$ of
index $i$. The differential $d:C^\Nov(v)_i\to C^\Nov(v)_{i-1}$ counts (with
signs) the number of flow lines from each critical point of index $i$ to each
critical point of index $i-1$ (in the universal cover of $M$). In contrast to
real-valued Morse theory, the image in $S^1$ of a flow-line may `wrap
around the circle' any number of times, so the differentials in
$C^\Nov(v)$ are formal power series in $z$. 

Now the torsion $\tau(\phi_v)$ lies in the image of the
canonical map 
\begin{equation*}
1+\Z\pi_\xi[[z]]z \to K_1(\Z\pi_\xi((z))).
\end{equation*}
For any $A$, let 
\begin{equation*}
W_1(A,\xi)=\Image\left(~1+A_\xi[[z]]z \to
K_1(A_\xi((z)))~\right).
\end{equation*}
Pajitnov and Ranicki showed~\cite{PajRan00} that $W_1(A,\xi)$ is a summand of
$K_1(A_\xi((z)))$ and is naturally isomorphic to the image of
$1+A_\xi[[z]]z$ in $K_1(A_\xi[[z]])$. In the light of
corollary~\ref{power_series_corollary} we have:
\begin{corollary}
\qquad$\displaystyle{W_1(A,\xi)\cong \frac{1+A_\xi[[z]]z}{C}}$. \hspace*{\fill}\qed
\end{corollary}
Pajitnov defined~\cite{Paj00} a logarithm based on the standard formula
\begin{equation*}
\cy{L}(1+\theta)=\theta-\frac{\theta^2}{2}+\frac{\theta^3}{3} \cdots 
\end{equation*}
which sends each element of $W_1(\Z\pi,\xi)$ to a formal power series
with one term for each conjugacy class $\beta$ of $\pi_1(M)$ such that
$f(\beta)\geq0$. He showed that the coefficient of $\beta$ in
$\cy{L}(\tau(\phi_v))$ is the number of closed orbits of the flow, counted with
signs, in the class $\beta$. \bigskip

This paper is structured as follows: In section~\ref{section:def_ex}
 we discuss local homomorphisms and augmentations and we define
 $K_1(A)$ and $\tEnd_0(A)$. Section~\ref{section:study_C} concerns the
 group $C$ of commutators which appears in
 theorems~\ref{Endomorphism_class_group}
 and~\ref{identify_K1}. Section~\ref{section:proof} is devoted to 
 the proofs of theorems~\ref{Endomorphism_class_group} and~\ref{identify_K1}.
\section{Definitions and Examples}
\label{section:def_ex}
Rings will be assumed associative with multiplicative unit. The set $\{0\}$
 will be considered a ring, in which $1=0$, but will not be considered
 a field.
\subsection{Local Homomorphisms}
\begin{definition}
\label{define_local}
A ring homomorphism $f:B\to A$ will be
called {\em local} if every square matrix $\alpha$ with entries in $B$
has the following property:
\begin{equation*}
\mbox{If $f(\alpha)$ is invertible then $\alpha$ is invertible.} 
\end{equation*}
Equivalently, $f$ is local if it has the property
that if $\alpha:P\to P$ is an endomorphism of a finitely
generated~(f.g.) projective $B$-module 
and the induced map $1\otimes\alpha:A\otimes_B P \to A\otimes_B P$ is
invertible then $\alpha:P\to P$ is invertible.
\end{definition}
For example, for any ring $B$ the surjection $B\to B/\rad(B)$ is a
local homomorphism (e.g.~Lam~\cite[Prop4.8 and (7) on p57]{Lam01}).

Formal power series provide another example.
If $A$ is any associative ring let $A[x]$ denote the polynomial
ring in a central indeterminate $x$. The augmentation
$\epsilon:A[x]\to A$ given by $\epsilon|_A=\id_A$ and
$\epsilon(x)=0$ induces local homomorphisms
$\displaystyle{\frac{A[x]}{(x^n)}\to A}$ and a local homomorphism
$\widehat{\epsilon}:A[[x]]\to A$; see
lemma~\ref{local_homo_cat_is_complete} below. 

It is easy to check that a composite of two local homomorphisms is
again local and that if a composite
$C\xrightarrow{g}B\xrightarrow{f}A$ is local then $g$ is local.

Although neither $A$ nor $B$ are local rings in most of our examples, 
we consider homomorphisms between local rings in our first
two lemmas to prove that definition~\ref{define_local} is
consistent with terminology used in algebraic geometry (see for
example Hartshorne~\cite[p73]{Har77}) and by Cohn~\cite[p388]{Coh85}.
Indeed, lemmas~\ref{B_local_rad_agree_implies_f_local}
and~\ref{local_gives_rad_condition} together imply that a
homomorphism $f:B\to A$ between local rings is a local homomorphism if
and only if $f^{-1}(\rad(A))=\rad(B)$.
\begin{lemma}\label{B_local_rad_agree_implies_f_local}
Suppose $f:B\to A$ is a ring homomorphism. If $B$ is a local ring, $A\neq
0$ and $f^{-1}(\rad(A))=\rad(B)$ then $f$ is a local homomorphism.
\end{lemma}
\noindent The hypotheses of the lemma are not redundant. For example
if $f$ is the inclusion of $\Z$ in $\Q$ then
$f^{-1}(\rad(\Q))=0=\rad(\Z)$ but $f$ is not a local homomorphism. On
the other hand, if $\Z_{(p)}$ denotes the local ring obtained from
$\Z$ by making invertible all the integers not divisible by the prime $p$ and
$f$ is the inclusion of $\Z_{(p)}$ in $\Q$ then $f^{-1}(\rad(\Q))\neq
\rad(\Z_{(p)})$ and $f$ is not a local homomorphism.
The following proof was pointed out to me by P.Ara.
\begin{proof}[Proof of lemma~\ref{B_local_rad_agree_implies_f_local}]
Suppose first that $B$ is a
division ring so the homomorphism $f:B\to A$ is an injection (since
$A\neq0$). Every homomorphism $B^n\to B^n$ is either an isomorphism or has
non-zero kernel so every square matrix $\alpha$ with entries in $B$ is
either invertible or a zero-divisor. It follows that $\alpha$ is
invertible in $A$ if and only if $\alpha$ is invertible in $B$ and
hence $f$ is local.

Now if $B$ is a local ring consider the commutative diagram
\begin{equation*}
\xymatrix@=2ex{
B\ar[r]^f\ar[d] & A\ar[d] \\
{\frac{B}{\rad(B)}} \ar[r] & {\frac{A}{\rad(A)}}.}
\end{equation*}

Since $B/\rad(B)$ is a division ring, the lower horizontal arrow is a local
homomorphism by the argument above. The vertical arrows are local
homomorphisms so it follows that $f$ is a local homomorphism.
\end{proof}
\begin{lemma}\label{local_gives_rad_condition}
If $f:B\to A$ is a local homomorphism then the preimage $f^{-1}(\rad(A))$
is a subset of $\rad(B)$. If, in addition, $A$ is a local ring then
$f^{-1}(\rad(A))= \rad(B)$ and $B$ is a local ring.
\end{lemma}
\noindent We remark that not every local homomorphism $f$ has
$f^{-1}(\rad(A))= \rad(B)$. For example, if $k$ is a field, the
inclusion of the formal power series ring $k[[x]]$ in the polynomial
extension $k[[x]][y]$ is a local homomorphism but
\begin{equation*}
f^{-1}(\rad(k[[x]][y]))=f^{-1}(0)=0\neq xk[[x]]=\rad(k[[x]]).
\end{equation*}
\begin{proof}[Proof of lemma~\ref{local_gives_rad_condition}]
An element $x\in B$ lies in $\rad(B)$ if and only if $1+bxb'$ is
invertible for all $b,b'\in B$ (e.g.~Lam~\cite[Lemma 4.3]{Lam01}).
If $x\in f^{-1}(\rad(A))$ then $f(x)\in\rad(A)$ so
$f(1+bxb')=1+f(b)f(x)f(b')$ is invertible. Since
$f$ is local it follows that $1+bxb'$ is invertible for all $b,b'\in
B$ and hence that $x\in\rad(B)$. Thus $f^{-1}(\rad(A))\subset
\rad(B)$. 

If $A$ is a local ring and $f$ is a local homomorphism then $x\in B$
is invertible if and only if $x\notin f^{-1}(\rad(A))$ so $B$ is a
local ring and $f^{-1}(\rad(A))=\rad(B)$.
\end{proof}
The next lemma says that any ring homomorphism can be made local in a
universal way:
\begin{lemma}
\label{localize_morphism}
Suppose $f:B\to A$ is a ring homomorphism. There is an initial object
\begin{equation}\label{make_homo_local}
B\xrightarrow{i_\Sigma} \Sigma^{-1}B\xrightarrow{f_\Sigma} A
\end{equation}
in the category of diagrams 
$B\xrightarrow{i} B'\xrightarrow{f'} A$
such that $f'i=f$ and $f'$ is local.
\end{lemma}
In other words, if $f:B\to A$ is the composite
$B\xrightarrow{i} B'\xrightarrow{f'} A$ and $f'$ is local then
there is a unique commutative diagram
\begin{equation}
\label{local_UP}
\begin{aligned}
\xymatrix @R=0.1ex{
& \Sigma^{-1}B \ar[dr]^{f_\Sigma}\ar[dd]_\gamma & \\
B\ar[ur]^{i_\Sigma} \ar[dr]_{i} &  & A \\
   & B' \ar[ur]_{f'} &
}
\end{aligned}
\end{equation}

Since it is an initial object, (\ref{make_homo_local}) is certainly unique.
The existence of (\ref{make_homo_local}) follows from the universal
localization of rings:
Given any set $\Sigma$ of matrices with entries in $B$, one can
adjoin formal inverses to each matrix in $\Sigma$,
to obtain a map $i_\Sigma:B\to \Sigma^{-1}B$ (see
Cohn~\cite[Ch.7]{Coh71} and Schofield~\cite[Ch.4]{Scho85}). Now $i_\Sigma$ is
$\Sigma$-inverting in the sense that, for every matrix
$\sigma\in\Sigma$, the image $i_\Sigma(\sigma)$ is
invertible. Moreover $i_\Sigma$ is characterized as the initial object
in the category of $\Sigma$-inverting homomorphisms $B\to B'$. In
other words every $\Sigma$-inverting homomorphism $B\to B'$ factors
through $i_\Sigma$ in a unique way. We remark that $i_\Sigma$ is an
epimorphism, i.e.~$fi_\Sigma=gi_\Sigma$ implies $f=g$.
\begin{proof}[Proof of lemma~\ref{localize_morphism}]
Given a ring homomorphism $f:B\to A$ let $\Sigma$ be
the set of $A$-invertible matrices with entries in $B$. In other
words, let $\sigma\in\Sigma$ if and only if $f(\sigma)$ is invertible.
Now $f$ is $\Sigma$-inverting, and hence factors uniquely through
$\Sigma^{-1}B$, i.e.~$f=f_\Sigma i_\Sigma$ where 
\begin{equation}
\label{make_homo_local_again}
B\xrightarrow{i_\Sigma} \Sigma^{-1}B \xrightarrow{f_\Sigma} A.
\end{equation}
It is proved in lemma 3.1 of~\cite{She01} that $f_\Sigma$ is a local
homomorphism.

To see that~(\ref{make_homo_local_again}) has the universal property
illustrated in~(\ref{local_UP}), suppose $f=f'i$ where $f':B'\to A$ is
local and $i:B\to B'$. If $\sigma\in\Sigma$ then $f(\sigma)$ is
invertible and, because $f'$ is
local, $i(\sigma)$ is invertible. It follows that there is a unique
map $\gamma:\Sigma^{-1}B\to B'$ such that $i=\gamma i_\Sigma$. Now
$f_\Sigma i_\Sigma=f'i=f'\gamma i_\Sigma$ and $i_\Sigma$ is an epimorphism
so we also have $f_\Sigma=f'\gamma$.%\hspace*{\fill}\qed
\end{proof}
One can obtain further examples of local homomorphisms by
limit constructions. For example a product of local
homomorphisms or an inverse limit of local homomorphisms is again a
local homomorphism. To make a general statement,
lemma~\ref{local_homo_cat_is_complete} below, let us briefly recall
the notion of limit in category theory. Suppose $F:J\to \cy{C}$ 
is a functor from a small category $J$ to a category $\cy{C}$. If $M$
is an object of $\cy{C}$ let $c_M:J\to \cy{C}$ denote the constant
functor which sends every object of $J$ to $M$ and every morphism to
$\id_M$. By definition a limit of $F$ is a final object in the
category of pairs $(M,\theta)$ where $M$ is an object of $\cy{C}$ and
$\theta$ is a natural transformation from $c_M$ to $F$.
The category $\cy{C}$ is said to be {\em complete} if every functor
$F:J\to\cy{C}$, where $J$ is small, has a limit.

The category of rings, for example, is complete (recall that we consider
$\{0\}$ a ring). 
Of interest here is the category in which an object is a ring
homomorphism and a morphism from $f:B\to A$ to
$f':B'\to A'$ is a commutative square:
\begin{equation*}
\xymatrix @=1.6ex{
B \ar[r]^{f}\ar[d] & A\ar[d] \\
B'\ar[r]_{f'} & A'
}.
\end{equation*}
This homomorphism category is also complete. Although the category of
local rings is not complete, one has
\begin{lemma}
\label{local_homo_cat_is_complete}
The category of local homomorphisms is complete.\hspace*{\fill}\qed
\end{lemma}
The proof of lemma~\ref{local_homo_cat_is_complete} is not difficult;
it suffices to check that equalizers and arbitrary products of local
homomorphisms are again local. The details are left to the reader.

Dually, one can attempt to construct examples of local homomorphisms by
colimit constructions. However, the category of local homomorphisms is
{\em not} cocomplete. 
For example, the coproduct of two copies of the local homomorphism
$\Z[x]/(x^2) \to \Z;x\mapsto0$  
is not local. On the positive side, the reader can check that a direct
colimit (often called a direct limit) of local homomorphisms is a
local homomorphism: 
\begin{lemma}
\label{direct_limits}
Let $J$ be a directed set. If $\left(\{f_j:B_j\to A_j\}_{j\in J}\;
,\; \{\theta^k_j:f_j\to f_k\}_{j\leq k}\right)$ is a 
direct system of local homomorphisms then the colimit
\begin{equation*}
\varinjlim f_j : \varinjlim B_j\to \varinjlim A_j
\end{equation*}
is a local homomorphism.\hspace*{\fill}\qed
\end{lemma}
\subsection{Local Augmentations}
The word `augmentation' is synonymous with `split surjection' or `retraction':
\begin{definition}
\label{define_augmentation}%
An {\em augmentation} $(\epsilon,j)$ is a pair of ring homomorphisms
\begin{equation*}
\epsilon:B\to A \quad \mbox{and}\quad j:A\to B
\end{equation*}
such that $\epsilon j=\id_A$. 
\end{definition}
The equation $\epsilon j=\id_A$ implies that $B$ can be expressed as a
direct sum $j(A)\oplus \Ker(\epsilon)$ of $(A,A)$-bimodules.
We shall usually suppress $j$, regarding $A$ as a subset of
$B$. Note that the category of augmentations is both complete and cocomplete.
In particular a direct or inverse limit of augmentations is again an
augmentation. A ring homomorphism which is both local and an
augmentation will be called a local augmentation.
\begin{lemma}
\label{nilpotent_implies_local}
If $\epsilon:B\to A$ is an augmentation and $I=\Ker(\epsilon)$
satisfies $I^n=0$ for some $n\in\N$ then $\epsilon$ is a local augmentation.
\end{lemma}
\begin{proof}
Suppose $\alpha$ is a square matrix with entries in $B$ and 
$\epsilon(\alpha)$ is invertible. The matrix 
$\alpha_0 = \epsilon(\alpha)^{-1}\alpha-1$ has entries in $I$ 
so $1+\alpha_0$ has inverse 
\begin{equation*}
1-\alpha_0 +\alpha_0^2- \cdots
+(-1)^{n-1}\alpha_0^{n-1}. 
\end{equation*}
Thus $\alpha$ is invertible and
\begin{equation*}
\alpha^{-1}=(1-\alpha_0 +\alpha_0^2- \cdots
+(-1)^{n-1}\alpha_0^{n-1})\epsilon(\alpha)^{-1}.%\tag*{\qed}
\qedhere
\end{equation*}
\end{proof}
\begin{lemma}[Localization and completion]
\label{loc_and_complete}
Suppose $\epsilon:B\to A$ is an augmentation with $\epsilon j=\id_A$. \begin{itemize}
\item[{\rm a)}] The universal localization $\epsilon_\Sigma:\Sigma^{-1}B\to
A$ is a local augmentation. %\\[0.3ex]\noindent
\item[{\rm b)}] Let $I=\Ker(\epsilon)$ and let
$\widehat{B}=\displaystyle{\varprojlim \frac{B}{I^n}}$ denote the
$I$-adic completion of $B$. For each $n$, the induced map
$\displaystyle{\epsilon_n:\frac{B}{I^n}\to A}$ is a local augmentation. The
inverse limit $\widehat\epsilon:\widehat{B}\to A$ is also a local
augmentation. %\\[0.3ex]\noindent
\item[{\rm c)}] 
There are natural isomorphisms 
\begin{equation*}
\Sigma^{-1}(\Sigma^{-1}B)\cong\Sigma^{-1}B;\quad
\Sigma^{-1}\widehat{B}\cong\widehat{B};\quad
\widehat{\widehat{B}}\cong\widehat{B};\quad
\widehat{\Sigma^{-1}B}\cong\widehat{B}.
\end{equation*} %\\[0.3ex]\noindent
\item[{\rm d)}] There is a natural homomorphism
$\gamma:\Sigma^{-1}B\to\widehat{B}$ such that
$\widehat{\epsilon}\gamma=\epsilon_\Sigma$.
\end{itemize}
\end{lemma}
\begin{proof}
a) By lemma 3.1 of~\cite{She01} we need only check that
$\epsilon_\Sigma$ is an augmentation. Indeed, 
$i_\Sigma j: A\to \Sigma^{-1}B$ has the property $\epsilon_\Sigma
i_\Sigma j =\epsilon j = \id_A$. \\[0.3ex] \noindent
b) Since an inverse limit of local augmentations is a local
augmentation, it suffices to show that $\epsilon_n:\frac{B}{I^n}\to A$ is
a local augmentation. The composite $A\xrightarrow{j}B\to
\frac{B}{I^n} \xrightarrow{\epsilon_n} A$ is
the identity morphism so $\epsilon_n$ is an augmentation. By
lemma~\ref{nilpotent_implies_local} $\epsilon_n$ is a local
homomorphism.  \\[0.3ex] \noindent
c) The first two identities merely reassert that $\epsilon_\Sigma$ and
$\widehat{\epsilon}$ are local. Now let
$\widehat{I}=\Ker(\widehat{\epsilon}:\widehat{B}\to A)$ and let
$I_\Sigma=\Ker(\epsilon_\Sigma:\Sigma^{-1}B\to A)$. The last two
identities follow respectively from natural isomorphisms
$\displaystyle{\frac{\widehat{B}}{{\widehat{I}}^{\raisebox{0.5ex}{$\scriptstyle{n}$}}}\cong \frac{B}{I^n}}$
and $\displaystyle{\frac{\Sigma^{-1}B}{I_\Sigma^{\raisebox{0.3ex}{$\scriptstyle{n}$}}}\cong
\frac{B}{I^n}}$. \\[0.3ex] \noindent
d) There are natural maps $\Sigma^{-1}B\to
\widehat{\Sigma^{-1}B}\to\widehat{B}$. %\hspace*{\fill}\qed
\end{proof}
\begin{example}
\label{twisted_power_series_example}
Let $X$ be any set and let $X^*$ denote the free monoid of words in
the alphabet $X$. Suppose $\xi:X\to \Aut(A);x\mapsto \xi_x$ is
any function which takes values in the group of ring automorphisms of $A$.
Let $A_\xi\langle X\rangle=A_\xi[X^*]$ denote the twisted monoid ring
of finite formal sums $\sum_{w\in X^*}a_w w$ with each $a_w\in
A$. Multiplication is defined by concatenation of words and 
by $ax=x\xi_x(a)$ for $a\in A$ and $x\in X$. 
Let $\epsilon:A_\xi\langle X\rangle\to A$ be the augmentation which
sends every letter in $X$ to zero.
By lemma~\ref{loc_and_complete}b) there is a local augmentation
$\widehat{\epsilon}:A_\xi\langle\langle X\rangle\rangle\to A$ where
$A_\xi\langle\langle X\rangle\rangle$ 
is the formal power series ring whose elements are infinite (formal)
sums $\sum_{w\in X^*}a_ww$. More generally, if $H\vartriangleleft
A_\xi\langle X\rangle$ is a (two-sided) ideal and $H\subset I=\Ker(\epsilon)$
then there is a local augmentation 
\begin{equation*}
\varprojlim_n \frac{A_\xi\langle
X\rangle}{I^n+H} \to A.
\end{equation*}
\end{example}
\subsection{$\tEnd_0(A)$ and $K_1(A)$}
\label{section:main_result}
Let $A$ be any associative ring. Let $\underline{\End}(A)$ denote the category
whose objects are pairs $(P,\alpha)$ where $P$ is a finitely generated
(f.g.) projective left $A$-module and $\alpha:P\to P$ is an
endomorphism. A morphism $(P,\alpha)\to (Q,\beta)$ in
$\underline{\End}(A)$ is a map $f:P\to Q$ such that $\beta f=f\alpha$.
\begin{definition}
\label{end_class_gp_defn}
Let $\tEnd_0(A)$ be the abelian group with one generator $[P,\alpha]$
for each isomorphism class of objects $(P,\alpha)$ in
$\underline{\End}(A)$ and relations
\begin{enumerate}
\item $[P',\alpha']=[P,\alpha]+[P'',\alpha'']$ if there is an exact
sequence 
\begin{equation*}
0\to(P,\alpha)\to(P',\alpha')\to(P'',\alpha'')\to0
\end{equation*}
in $\underline{\End}(A)$.
\item $[P,0]=0$ for all f.g.~projective modules $P$.
\end{enumerate}
\end{definition}
Definition~\ref{end_class_gp_defn} is consistent with
definition~\ref{end_class_gp_free}; see~\cite[Lemma 2.4]{She01}. To
define $K_1(A)$, let $\underline{\Aut}(A)$ denote the full subcategory
of $\underline{\End}(A)$ whose objects are pairs $(P,\alpha)$ such
that $\alpha$ is an automorphism.
\begin{definition}
\label{define_K1}
The abelian group $K_1(A)$ is generated by the
isomorphism classes $[P,\alpha]$ in $\underline{\Aut}(A)$ subject to
relations:
\begin{enumerate}
\item $[P\oplus P',\alpha\oplus \beta]=[P,\alpha]+[P,\beta]$.
\item\label{product_relation} $[P,\alpha\beta]=[P,\alpha]+[P,\beta]$.
\end{enumerate}
\end{definition}
The group $K_1(A)$ is also  unchanged if one replaces the projective
modules in definition~\ref{define_K1} by free modules throughout. Indeed,
relation~\ref{product_relation}.~implies that
$[P,1_P]=[P,1_P]+[P,1_P]$ and hence $[P,1_P]=0\in K_1(A)$
for all $P$. It follows that $[P,\alpha]$ can be identified with $[P\oplus
Q,\alpha\oplus 1_Q]$ where $P\oplus Q$ is free.

There is a natural isomorphism between $K_1(A)$ and
$\displaystyle{\GL(A)^\ab}$, the abelianization of the direct limit
$\GL(A)=\varinjlim \GL_n(A)$ of general linear groups over $A$
(e.g.\cite[p109]{Sil81},\cite[Thm 3.1.7]{Ros94}).
The commutator subgroup $[\GL(A),\GL(A)]$ is generated by the `elementary'
matrices $b_{ij}(a)=[b_{ik}(a),b_{kj}(1)]$ where $b_{ij}(a)$ is the
matrix whose $ij$th entry is $a$, whose diagonal entries are $1$ and
whose other entries are $0$.
\section{The commutator group $C$}
\label{section:study_C}%
Suppose $\epsilon:B\to A$ is a local augmentation. In this section we
study the group 
\begin{equation*}
C=\langle(1+ab)(1+ba)^{-1}\mid
\epsilon(ab)=\epsilon(ba)=0\rangle\subset\epsilon^{-1}(1).
\end{equation*}
A consequence of proposition~\ref{comparing_commutators} below is that
$[\epsilon^{-1}(1),\epsilon^{-1}(1)]\subset C$, a fact we shall need
in the proof of theorem~\ref{identify_K1}. 
Although $[\epsilon^{-1}(1),\epsilon^{-1}(1)]\neq C$ in general,
the following lemma says that the image of $C$ in $K_1(B)$ is trivial.
\begin{lemma}\label{C_dies_in_K1}
Suppose $B$ is any ring and $a,b\in B$. The element $1+ab$ is
invertible if and only if $1+ba$ is invertible. Moreover, if $1+ab$
and $1+ba$ are invertible then $[1+ab]=[1+ba]\in K_1(B)$.
\end{lemma}
\begin{proof}% One can check that
$\displaystyle{\left(\begin{matrix}
1 &           -a \\
0 & \phantom{-}1
\end{matrix}\right)
\left(\begin{matrix}
1+ab & 0 \\
0 & 1
\end{matrix}\right)
\left(\begin{matrix}
1     & 0 \\
b & 1
\end{matrix}\right)
=
\left(\begin{matrix}
1     & 0 \\
b & 1
\end{matrix}\right)
\left(\begin{matrix}
1 & 0 \\
0 & 1+ba
\end{matrix}\right)
\left(\begin{matrix}
1 & -a \\
0 & \phantom{-}1
\end{matrix}\right)}$.
%\hspace*{\fill}\qed
\end{proof}
\begin{notation}
Let $B^\bullet$ and $A^\bullet$ denote the
groups of units in $B$ and $A$ respectively. Suppose $S$ is a subset
of $B\times B$ with the property that $1+ba\in B^\bullet$ for all $(a,b)\in S$. 
Let $C(S)$ denote the intersection
of $\epsilon^{-1}(1)$ with the group generated by 
$\{(1+ab)(1+ba)^{-1}\mid (a,b)\in S\}$. In symbols
\begin{equation*}
C(S)=\epsilon^{-1}(1)\cap\langle(1+ab)(1+ba)^{-1}\mid (a,b)\in S\rangle.
\end{equation*}
\end{notation}
We usually describe $S$ in terms of
equations or conditions on $a$ and $b$.
For example, the commutator group in theorem~\ref{identify_K1} is
$C=C(\epsilon(ab)=\epsilon(ba)=0)$.
Note that if $S^T=\{(a,b)\mid(b,a)\in S\}$ then
$C(S)=C(S^T)=C(S\cup S^T)$
because $(1+ab)(1+ba)^{-1}=\left((1+ba)(1+ab)^{-1}\right)^{-1}$.
The following fact was attributed to L.Vaserstein by
V.Srinivas~\cite[p5]{Sri96}.
\begin{lemma}
\label{adding_commutation}
If $1+ac$ and $1+ba$ are invertible and $ac=ca$ then
\begin{equation*}
(1+ab)(1+ba)^{-1}=(1+a(b+c+bac))(1+(b+c+bac)a)^{-1}.
\end{equation*}
\end{lemma}
\begin{proof}
Observe that $(1+ab)(1+ac)=1+a(b+c+bac)$ and
$(1+ba)(1+ca)=1+(b+c+bac)a$.
%\hspace*{\fill}\qed
\end{proof}
\begin{proposition}\label{comparing_commutators}
Let $\epsilon:B\to A$ be a local augmentation.
\begin{enumerate}
\item $C(\epsilon(a)=0)~=~C(\epsilon(ab)=\epsilon(ba)=0)~=~C(\epsilon(ba)=0)~=~C(1+ba\in
B^\bullet)$.\label{commut_all}%
\item $[B^\bullet,B^\bullet]\cap \epsilon^{-1}(1)~=~C(b\in B^\bullet;1+ba\in B^\bullet)$\label{commut_b_1+ba}
\item
$[\epsilon^{-1}(1),\epsilon^{-1}(1)]~=~C(\epsilon(a)=0;\epsilon(b)=\zeta)$
for any $\zeta\in A$ which commutes with every element of
$I=\Ker(\epsilon)$.\label{commut_0_0}%\\[1ex]
\end{enumerate}
In particular,
\begin{equation}
\label{commutator_containment}%
[\epsilon^{-1}(1),\epsilon^{-1}(1)]~\subset~
[B^\bullet,B^\bullet]\cap\epsilon^{-1}(1)~\subset~C~=~C(\epsilon(ab)=\epsilon(ba)=0)
\end{equation}
and $C$ is a normal subgroup of $\epsilon^{-1}(1)$ with abelian quotient.
\end{proposition}
Statement~\ref{commut_all}.~will follow from the proof of
theorem~\ref{identify_K1} in section~\ref{section:proof}
below. We do not use \ref{commut_all}.~in the proof of
\ref{commut_b_1+ba}.~or \ref{commut_0_0}. 
%
%*
\begin{proof}[Proof of statement \ref{commut_b_1+ba}.]
(Compare Silvester~\cite[p135]{Sil81}) Suppose $\alpha,\beta\in
B^\bullet$.
Define $a=\alpha\beta - \alpha$ and $b=\alpha^{-1}$ so that
$(1+ab)(1+ba)^{-1}=\alpha\beta\alpha^{-1}\beta^{-1}$.
Now $b\in B^\bullet$ and $1+ba=\beta\in B^\bullet$ so
$[B^\bullet,B^\bullet]\cap\epsilon^{-1}(1)\subset 
C(b\in B^\bullet;1+ba\in B^\bullet)$. 

Inversely, if $a$,$b\in B$ with $b\in B^\bullet$ and $1+ba\in
B^\bullet$ one can set $\alpha=b^{-1}$ and $\beta=1+ba$ 
to obtain 
$\alpha\beta\alpha^{-1}\beta^{-1}=(1+ab)(1+ba)^{-1}$.
%\hspace*{\fill}\qed
\end{proof}
%*
\begin{proof}[Proof of statement \ref{commut_0_0}.]
We continue to use the notation of the preceding proof. Observe that
$\epsilon(\alpha)=\epsilon(\beta)=1$ if and only if $\epsilon(a)=0$
and $\epsilon(b)=1$. Thus
$[\epsilon^{-1}(1),\epsilon^{-1}(1)]=C(\epsilon(a)=0;\epsilon(b)=1)$. 
It remains to prove that
\begin{equation*}
C(\epsilon(a)=0;\epsilon(b)=1)~=~C(\epsilon(a)=0;\epsilon(b)=\zeta).
\end{equation*} 
If $\epsilon(a)=0$ and $\epsilon(b)=1$, put $c=\zeta-1$ so that
$\epsilon(b+c+bac)=\zeta$. Lemma~\ref{adding_commutation} implies that
$C(\epsilon(a)=0;\epsilon(b)=1)~\subset~C(\epsilon(a)=0;\epsilon(b)=\zeta)$.
Conversely, if $\epsilon(a)=0$ and $\epsilon(b)=\zeta$, put
$c=1-\zeta$ so that $\epsilon(b+c+bac)=1$ and
lemma~\ref{adding_commutation} implies that
$C(\epsilon(a)=0;\epsilon(b)=\zeta)~\subset~C(\epsilon(a)=0;\epsilon(b)=1)$.
%\hspace*{\fill}\qed 
\end{proof}
Using lemma~\ref{adding_commutation} one can show that the second
inclusion of (\ref{commutator_containment}) is an equality in certain cases of interest (e.g.~if $A$ is
a local ring). However, neither inclusion is an equality in general as
the following example illustrates.
\begin{example}\label{inclusions_are_proper}
Suppose $B=A[[x]]$, where $x$ is a central indeterminate, and
$\epsilon$ is the local augmentation $A[[x]]\to A;x\mapsto 0$. Pajitnov and
Ranicki observed~\cite[\S3.2]{PajRan00} that the canonical map 
$\epsilon^{-1}(1)^\ab\to K_1(B)$ need not be injective. In fact,
whenever there exist elements $a,b\in A$ such that $ab\neq ba$ one
finds that
\begin{equation*}
(1+a(bx))(1+(bx)a)^{-1}=1+(ab-ba)x + \cdots
\end{equation*}
is an element of $C$ but is not in $[\epsilon^{-1}(1),\epsilon^{-1}(1)]$.
Furthermore, if one can choose $a\in A^\bullet$ and $b\in
B$ with $ab\neq ba$ then 
\begin{equation*}
a(1+bx)a^{-1}(1+bx)^{-1}=1+(aba^{-1}-b)x + \cdots
\end{equation*}
lies in $[B^\bullet,B^\bullet]\cap\epsilon^{-1}(1)$ but not in
$[\epsilon^{-1}(1),\epsilon^{-1}(1)]$. On the other hand, if say
$A=\Z\langle y,z\rangle$ is the free associative ring on two
generators then $A^\bullet=\{\pm1\}$ so 
$B^\bullet=\pm(1+A[[x]]x)$ and $[B^\bullet,B^\bullet]\subset
1+A[[x]]x^2$. Putting $a=xz$ and $b=y$ we have 
\begin{equation*}
(1+ab)(1+ba)^{-1}=1+(zy-yz)x + \cdots
\end{equation*}
so $(1+ab)(1+ba)^{-1}$ lies in $C$ but not in $[B^\bullet,B^\bullet]$.
\end{example}
It follows from lemma~\ref{loc_and_complete}d),
proposition~\ref{comparing_commutators} and
example~\ref{inclusions_are_proper} that, in the context of
theorem~\ref{Endomorphism_class_group}, we have
$[\epsilon_\Sigma^{-1}(1),\epsilon_\Sigma^{-1}(1)]
\subsetneqq C$ in $\Sigma^{-1}A[x]$ whenever $A$ is non-commutative.
\section{Proof of Theorems~\ref{Endomorphism_class_group}
and~\ref{identify_K1}}
\label{section:proof}
Suppose $\epsilon:B\to A$ is a local homomorphism and $j:A\to B$
satisfies $\epsilon j=\id_A$. 
We use the same symbols $\epsilon$ and $j$ to denote the functors
$A\otimes_B \functor$ and $B\otimes_A \functor$ and the induced maps
$K_1(B)\to K_1(A)$ and $K_1(A)\to K_1(B)$.
Since $\epsilon j= \id_A$ induces the identity on $K_1(A)$ 
we have a decomposition
\begin{equation*}
K_1(B)=K_1(A)\oplus \tK_1(B)
\end{equation*}
where, by definition, $\tK_1(B)=\Ker(\epsilon:K_1(B)\to K_1(A))$. To
prove theorem~\ref{identify_K1} 
we must show that $\tK_1(B)$ is isomorphic to
$\displaystyle{\epsilon^{-1}(1)/C}$ where $C$ is the subgroup
of $\epsilon^{-1}(1)$ generated by
the subset  
\begin{equation}
\left\{(1+ab)(1+ba)^{-1} \bigm{|}
a,b \in B,
\epsilon(ab)=\epsilon(ba)=0\right\}. \label{comm_rel}
\end{equation}
We shall continue to write the group operation multiplicatively
in $\epsilon^{-1}(1)/C$ but additively in $K_1(B)$.

We first deduce theorem~\ref{Endomorphism_class_group} from
theorem~\ref{identify_K1}:
%
%*
\begin{proof}[Proof of theorem~\ref{Endomorphism_class_group}.]
Recall that $\epsilon:A[x]\to A;x\mapsto0$ and $\Sigma$ denotes the
set of matrices $\sigma$ with entries in $A[x]$ such that
$\epsilon(\sigma)$ is invertible. The universal localization
$\epsilon_\Sigma:\Sigma^{-1}A[x]\to A$ is a local augmentation by
lemma~\ref{loc_and_complete}a). Ranicki showed~\cite[Prop10.21]{Ran98}
that there is an isomorphism
\begin{equation*}
\displaystyle{\tEnd_0(A)\cong\tK_1(\Sigma^{-1}A[x]);[P,\alpha]\mapsto[1-\alpha
x]}.
\end{equation*}
By theorem~\ref{identify_K1}, the map
$\displaystyle{\tEnd_0(A) \to \epsilon_\Sigma^{-1}(1)/C; 
[P,\alpha] \mapsto D(1-\alpha x)}$
is also an isomorphism. 
%\hspace*{\fill}\qed
\end{proof}
\begin{proof}[Outline proof of theorem~\ref{identify_K1}]
The proof of theorem~\ref{identify_K1} is analogous to the proof of
the identity $K_1(k)\cong k^\bullet$, where $k$ is a (commutative)
field. The latter is proved by observing that the determinant   
\begin{align*}
\GL(k) &\to k^\bullet \\
\alpha &\mapsto \det(\alpha) 
\end{align*}
induces a map $\det:K_1(k)\to k^\bullet$ which is inverse to the
canonical map from $k^\bullet$ to $K_1(k)$. 

In a non-commutative setting this traditional determinant is not defined.
By lemma~\ref{C_dies_in_K1} there is a canonical map
$\epsilon^{-1}(1)/C\to \tK_1(B)$; we shall construct an inverse
$D:\tK_1(B)\to\epsilon^{-1}(1)/C$ which is a 
version of the Dieudonn\'e determinant. The idea is
that the class in $K_1(B)$ represented by a given invertible matrix is
unchanged when one performs elementary operations, adding a (left)
multiple of one row to another row, or adding a (right) multiple of one column 
to another column. If $[\alpha]\in\tK_1(B)$ then $\alpha$ can be
 reduced to a diagonal matrix by these row and column operations.
One can then define $D(\alpha)$ to be the product of the diagonal
 entries in this diagonal matrix. Modulo appropriate relations, which
turn out to be $C$, the determinant $D$ is well-defined.
%\hspace*{\fill}\qed
\end{proof}
We treat the preceding proof in more detail below showing in fact that $D$ is
well-defined modulo $C_0~=~C(\epsilon(a)=0)~=~\langle(1+ab)(1+ba)^{-1}\mid
\epsilon(a)=0\rangle$. Let us now deduce that the various definitions
of $C$ in part~1.~of proposition~\ref{comparing_commutators} are identical:
%*
\begin{proof}[Proof of part~\ref{commut_all}.~of proposition~\ref{comparing_commutators}]
It is immediate that
\begin{equation}
\label{commut_all_inclusions}
\begin{split}
C_0~=~C(\epsilon(a)=0)~\subset~C~=C&(\epsilon(ab)=\epsilon(ba)=0) \\
&\subset~C(\epsilon(ba)=0)~\subset~C(1+ba\in
B^\bullet).
\end{split}
\end{equation}
By lemma~\ref{C_dies_in_K1}, $C(1+ba\in B^\bullet)$ vanishes in
$\tK_1(B)$ so the composite
\begin{equation*}
\frac{\epsilon^{-1}(1)}{C_0}\twoheadrightarrow\frac{\epsilon^{-1}(1)}{C(1+ba\in
B^\bullet)}\to\tK_1(B)\to\frac{\epsilon^{-1}(1)}{C_0}
\end{equation*}
is the identity. Hence $C_0=C(1+ba\in
B^\bullet)$ and all the inclusions in (\ref{commut_all_inclusions})
are equalities. In particular, $C_0=C$.
%\hspace*{\fill}\qed
\end{proof}
We begin now a detailed proof
that $D:\tK_1(B)\to \epsilon^{-1}(1)/C_0$ is
well-defined. Let $\tAuto(B)$ denote the full subcategory of
$\underline{\End}(B)$ whose objects are pairs $(P,\alpha)$ such that
$\alpha:P\to P$ is an automorphism and $\epsilon(\alpha)=1$.
\begin{lemma}
\label{genrel_tK}
$\tK_1(B)$ is (isomorphic to) the abelian group generated by
isomorphism classes $[P,\alpha]$ of automorphisms
$(P,\alpha)\in\tAuto(B)$ subject to the following relations:
\begin{enumerate}
\item $[P\oplus P',\alpha\oplus \beta]=[P,\alpha]+[P',\beta]$
\item $[P,\alpha\beta]=[P,\alpha]+[P,\beta]$
\item \label{extra_relation} If $(P,c)=j(P_0,c_0)\in\underline{\Aut}(B)$ and 
$(P,\alpha)\in\tAuto(B)$ then $[P,c\alpha c^{-1}]=[P,\alpha]$.
\end{enumerate}
\end{lemma}
In lemma~\ref{genrel_tK}, just as in
definition~\ref{define_K1}, one may replace f.g.~projective modules by
f.g.~free modules throughout. We usually abbreviate $[P,\alpha]$
to $[\alpha]$.
%*
\begin{proof}[Proof of lemma~\ref{genrel_tK}]
Let $T$ denote the abelian group with the generators and relations
given in lemma~\ref{genrel_tK}. There is an obvious map $T\to K_1(B)$
and the composite $T\to K_1(B)\xrightarrow{\epsilon} K_1(A)$ is zero,
so the image of $T$ lies in $\tK_1(B)$. Conversely, given a generator
$[\alpha]\in K_1(B)$ we may write $\alpha=c\talpha$ where
$c=j\epsilon(\alpha)$ and $\talpha=c^{-1}\alpha$ has
the property $\epsilon(\talpha)=1$. We can define a map 
\begin{equation*}
\theta:K_1(B)\to T;\quad [\alpha]\mapsto [\talpha]=[c^{-1}\alpha]
\end{equation*}
Plainly $\theta[\alpha\oplus \beta]=\theta[\alpha]+
\theta[\beta]$. We must check also that
$\theta[\alpha\beta]=\theta[\alpha]+\theta[\beta]$. Indeed, if
$c=\epsilon(\alpha)$ and $d=\epsilon(\beta)$ then
\begin{equation*}
\theta[\alpha\beta]=\theta[c\talpha d\tbeta]=\theta[(cd)d^{-1}\talpha
d\tbeta] = [d^{-1}\talpha
d\tbeta] = [d^{-1}\talpha
d]+[\tbeta]
\end{equation*}
so by relation~\ref{extra_relation}.~we have
$\theta[\alpha\beta]=
[\talpha]+[\tbeta]=\theta[\alpha]+\theta[\beta]$. Plainly the
restriction of $\theta$ to $\tK_1(B)$ is
inverse to the canonical map $T\to \tK_1(B)$. 
%\hspace*{\fill}\qed
\end{proof}

It is easy to define a map 
$i:\epsilon^{-1}(1)\to \tK_1(B)$. Indeed, a unit $\alpha\in B^\bullet$
determines an automorphism $B\to B;x\mapsto x\alpha$ of the free
$B$-module on one generator. Moreover, if $\epsilon(\alpha)=1$ then
$[\alpha]\in \tK_1(B)\subset K_1(B)$. 

We prove next that $i$ is surjective, the idea being to reduce an
automorphism of $B^n$, `by row and column operations', to a direct sum
of automorphisms of smaller modules. If $(B,\alpha)\in\tAuto(B)$
then $[\alpha]$ is certainly in the image of $i$.
\begin{lemma}
\label{LDUform}
If $\alpha$ is an automorphism of $B^n$ with $n\geq2$ and
$\epsilon(\alpha)=1$ then there is a unique expression
\begin{equation}
\label{LDU}
\alpha= 
\left(\begin{matrix}
1  &  0 \\
l  &  1
\end{matrix}\right)
\left(\begin{matrix}
d_1 & 0 \\
0   & d_2
\end{matrix}\right)
\left(\begin{matrix}
1 & u \\
0 & 1    
\end{matrix}\right)
\end{equation}
where $l:B\to B^{n-1}$,
$u:B^{n-1}\to B$, $d_1:B\to B$ and $d_2:B^{n-1}\to B^{n-1}$. This
expression has the properties
$\epsilon(d_1)=1_A$, $\epsilon(d_2)=1_{A^{n-1}}$, $\epsilon(u)=0$
and $\epsilon(l)=0$.
\end{lemma}
%
%*
\begin{proof}[Proof of lemma~\ref{LDUform}]
We may write $\alpha=\left(\begin{matrix}
\alpha_{11} & \alpha_{12} \\
\alpha_{21} & \alpha_{22}
\end{matrix}\right)$
where 
\begin{equation*}
\epsilon(\alpha_{11})=1,\quad \epsilon(\alpha_{22})=1,\quad
\epsilon(\alpha_{12})=0 \quad\mbox{and}\quad \epsilon(\alpha_{21})=0.
\end{equation*}
Now 
$\alpha_{11}$ is invertible, since $\epsilon$ is a local homomorphism,
so
\begin{equation*}
\alpha=\left(\begin{matrix}
1                             & 0 \\ 
\alpha_{21}\alpha_{11}^{-1}   & 1
\end{matrix}\right)
\left(\begin{matrix}
\alpha_{11}              & 0 \\
     0                   &
     \alpha_{22}-\alpha_{21}\alpha_{11}^{-1}\alpha_{12}
\end{matrix}\right)
\left(\begin{matrix}
1  & \alpha_{11}^{-1}\alpha_{12} \\
0  & 1
\end{matrix}\right).
\end{equation*}
This equation proves existence of the expression~(\ref{LDU}) and the properties
in the last sentence of the lemma follow immediately. To show
uniqueness suppose
\begin{equation*}
\left(\begin{matrix}
1  &  0 \\
l  &  1
\end{matrix}\right)
\left(\begin{matrix}
d_1 & 0 \\
0        & d_2
\end{matrix}\right)
\left(\begin{matrix}
1 & u \\
0 & 1    
\end{matrix}\right)
=
\left(\begin{matrix}
1         &  0 \\
l'  &  1
\end{matrix}\right)
\left(\begin{matrix}
d'_1 & 0 \\
0        & d'_2
\end{matrix}\right)
\left(\begin{matrix}
1 & u' \\
0 & 1    
\end{matrix}\right)~.
\end{equation*}
It follows that
\begin{equation*}
\left(\begin{matrix}
1                     &  0 \\
-l' + l &  1
\end{matrix}\right)
\left(\begin{matrix}
d_1 & 0 \\
0        & d_2
\end{matrix}\right)
=
\left(\begin{matrix}
d'_1 & 0 \\
0        & d'_2
\end{matrix}\right)
\left(\begin{matrix}
1 & u'-u \\
0 & 1    
\end{matrix}\right)
\end{equation*}
which implies first that $d_1=d'_1$ and $d_2=d'_2$
and second, since $d_1$ and $d_2$ are invertible, that $u=u'$ and
$l=l'$.
%\hspace*{\fill}\qed
\end{proof}
Since $\left[\begin{matrix}
1  &  0 \\
l  &  1
\end{matrix}\right]=
\left[\begin{matrix}
1  &  u \\
0  &  1
\end{matrix}\right]=0\in\tK_1(B)$ 
the `existence' part of lemma~\ref{LDUform} implies by induction that
the natural map $\epsilon^{-1}(1)\to \tK_1(B)$ is surjective.

The `uniqueness' part of lemma~\ref{LDUform} leads to a map from the
objects in $\tAuto(B)$ to $\epsilon^{-1}(1)/{C_0}$, defined
recursively. This map is analogous to the Dieudonn\'e
determinant~\cite{Die43} (compare also Klingenberg~\cite{Kli62},
Draxl~\cite[Ch20]{Dra83}, Silvester~\cite[pp122-140]{Sil81},
Srinivas~\cite[Example 1.6]{Sri96}).
We continue to use the notation of lemma~\ref{LDUform}.
\begin{definition}\label{Define_D}
Suppose $\alpha:B^n\to B^n$ is an automorphism. If $n=1$ then define
$D(\alpha)=\alpha \in \epsilon^{-1}(1)/C_0$. If $n\geq2$ then define
\begin{equation*}
D(\alpha)=d_1D(d_2)=\alpha_{11}D(\alpha_{22}-\alpha_{21}\alpha_{11}^{-1}\alpha_{12})
\in \frac{\epsilon^{-1}(1)}{C_0}.
\end{equation*}  
\end{definition} 
It is easy to see that if $\alpha=\id:B^n\to B^n$ then
$D(\alpha)=1$. One can extend the definition of $D$ to automorphisms of
finitely generated projective modules by $D(\alpha:P\to P)=D(\alpha\oplus
1_Q)$ where $P\oplus Q$ is finitely generated and free.
 
Now if we can show that $D$ respects the relations 1-3 of
lemma~\ref{genrel_tK} then $D$ will induce a map 
\begin{equation*}
\tK_1(B)\to
\epsilon^{-1}(1)/C_0
\end{equation*}
which is plainly inverse to
$i:\epsilon^{-1}(1)/C_0\to \tK_1(B)$.
It suffices to consider automorphisms of free modules; we shall check the
relations by induction starting with relation 1.: 
\begin{lemma}
If $\alpha:B^n\to B^n$ and $\beta:B^m\to B^m$ then
$D(\alpha\oplus\beta)=D(\alpha)D(\beta)$.
\end{lemma}
\begin{proof}
We perform induction on $n$. The case $n=1$ is trivial so let us suppose
$n\geq2$. Writing
$\alpha=\left(\begin{matrix}
\alpha_{11} & \alpha_{12} \\
\alpha_{21} & \alpha_{22}
\end{matrix}\right)$ with the $\alpha_{ij}$ defined as in the proof of
lemma~\ref{LDUform}, we have
\begin{align*}
D(\alpha\oplus\beta) &=D\left(\begin{matrix} 
\alpha_{11}& \alpha_{12} & 0 \\
\alpha_{21}& \alpha_{22} & 0 \\
0          & 0           & \beta
\end{matrix}\right) \\
&=\alpha_{11}
D\left(\left(\begin{matrix}\alpha_{22} & 0 \\
0 & \beta \end{matrix}\right)
-
\left(\begin{matrix}
\alpha_{21} \\
0
\end{matrix}\right)
\alpha_{11}^{-1}
\bigl(\begin{matrix}
\alpha_{12} & 0
\end{matrix}\bigr)
\right) \\
&=
\alpha_{11}D\left(\begin{matrix}
\alpha_{22}-\alpha_{21}\alpha_{11}^{-1}\alpha_{12} & 0 \\
0 & \beta
\end{matrix}\right) \\
&=
\alpha_{11}D(\alpha_{22}-\alpha_{21}\alpha_{11}^{-1}\alpha_{12})D(\beta)
\quad\mbox{(by the inductive hypothesis)}  \\
&=D(\alpha)D(\beta) 
%\tag*{\qed}
\qedhere
\end{align*}
\end{proof}
We treat relations 2.~and 3.~together:
\begin{proposition}
\label{D_respects_commutation}
\begin{itemize}
\item[a)] If $\alpha:B^n\to B^n$ and $\beta:B^n\to B^n$ are automorphisms and
$\epsilon(\alpha)=\epsilon(\beta)=1:A^n\to A^n$ then
$D(\alpha\beta)=D(\alpha)D(\beta) \in \epsilon^{-1}(1)/C_0$. 
\item[b)] Suppose $\alpha:B^n\to B^m$,
$\beta:B^m\to B^n$, and either $\epsilon(\alpha)=0$ or
$\epsilon(\beta)=0$ (or both).
Then
$D(1+\alpha\beta)=D(1+\beta\alpha)\in\epsilon^{-1}(1)/C_0$.
\end{itemize}
\end{proposition}
Part a) of proposition~\ref{D_respects_commutation} says that $D$
respects relation 2.~of lemma~\ref{genrel_tK}. Part b) implies that $D$
respects relation 3.~for, if $\alpha=1+\alpha_0$ and
$\epsilon(\alpha_0)=0$ then 
$D(c\alpha c^{-1})=D(1+c\alpha_0 c^{-1})=D(1+\alpha_0
c^{-1}c)=D(\alpha)$.
%
%*
\begin{proof}[Proof of proposition~\ref{D_respects_commutation}]
We shall prove parts a) and b) at the same time, by induction. In part
a) there is one statement, denoted $\mfa(n)$, for each positive integer
$n$. In part b) there is one statement $\mfb(m,n)$ for each pair
$(m,n)$ of positive integers and we have $\mfb(m,n)\Leftrightarrow \mfb(n,m)$.
 
Plainly $\mfa(1)$ and $\mfb(1,1)$ are true. We shall establish the following
inductive steps:
\begin{align}
\mfb(1,n-1) \wedge \mfa(n-1) &\Rightarrow \mfa(n)~\mbox{for all $n\geq2$} \label{a_induction}\\
 \mfa(m) \wedge \mfb(1,m) \wedge \mfb(m,n-1) &\Rightarrow \mfb(m,n)~\mbox{for
 $m\geq1$, $n\geq2$}\label{b_induction}
\end{align}
These are sufficient to prove $\mfa(n)$ and $\mfb(m,n)$ for
all $m$ and $n$. Indeed, a special case of (\ref{b_induction}) reads
$\mfa(1)\wedge \mfb(1,1)\wedge \mfb(1,n-1) \Rightarrow \mfb(1,n)$
so $\mfb(1,n)$ holds for all $n\geq1$. It then follows from
(\ref{a_induction}) that $\mfa(n)$ holds for all $n\geq1$. Finally
(\ref{b_induction}) shows, by induction on $n$, that $\mfb(m,n)$ holds
for all $m,n\geq1$.
\smallskip

\noindent {\em Proof of (\ref{a_induction}).}
Suppose $\alpha,\beta:B^n\to B^n$.
By lemma~\ref{LDUform} it suffices to
show $D(\alpha\beta)=D(\alpha)D(\beta)$ in the cases \\ \noindent
i) $\alpha=\left(\begin{matrix}
1  &  0 \\
l  &  1
\end{matrix}\right)$,\quad ii) $\alpha=\left(\begin{matrix}
d_1 & 0 \\
0        & d_2
\end{matrix}\right)$\quad and \quad iii) $\alpha=\left(\begin{matrix}
1 & u \\
0 & 1    
\end{matrix}\right)$. \\ \noindent Here $l:B\to B^{n-1}$, $d_1:B\to B$,
$d_2:B^{n-1}\to B^{n-1}$ and $u:B^{n-1}\to B$. \\[2mm]
\noindent
i) Let us write $\beta=\left(\begin{matrix}
1         &  0 \\
l'  &  1
\end{matrix}\right)
\left(\begin{matrix}
d'_1 & 0 \\
0        & d'_2
\end{matrix}\right)
\left(\begin{matrix}
1 & u' \\
0 & 1    
\end{matrix}\right)$.
Now
\begin{equation*}
D(\alpha\beta)=D\left(\begin{matrix}
1  &  0 \\
l+l'  &  1
\end{matrix}\right)
\left(\begin{matrix}
d'_1 & 0 \\
0        & d'_2
\end{matrix}\right)
\left(\begin{matrix}
1 & u' \\
0 & 1    
\end{matrix}\right)=d'_1D(d'_2)=D(\beta)=D(\alpha)D(\beta).
\end{equation*}
\vspace*{-2ex}
\begin{flalign*}
\mbox{ii)}\quad 
D(\alpha\beta) &=D\left(
\left(\begin{matrix}
d_1 & 0 \\
0        & d_2
\end{matrix}\right)
\left(\begin{matrix}
1  &  0 \\
l'  &  1
\end{matrix}\right)
\left(\begin{matrix}
d'_1 & 0 \\
0        & d'_2
\end{matrix}\right)
\left(\begin{matrix}
1 & u' \\
0 & 1    
\end{matrix}\right)
\right) && \\
&=
D\left(\left(\begin{matrix}
1  &  0 \\
d_2l'd_1^{-1}  &  1
\end{matrix}\right)
\left(\begin{matrix}
d_1d'_1 & 0 \\
0        & d_2d'_2
\end{matrix}\right)
\left(\begin{matrix}
1 & u' \\
0 & 1    
\end{matrix}\right)\right) && \\
&=d_1d'_1D(d_2d'_2).
\end{flalign*}
Now $\epsilon^{-1}(1)/C_0$ is abelian by
part~\ref{commut_0_0}. of proposition~\ref{comparing_commutators} so
by $\mfa(n-1)$ we have
$D(\alpha\beta)=d_1d'_1D(d_2)D(d'_2)=d_1D(d_2)d'_1D(d'_2)=D(\alpha)D(\beta)$.
This completes part ii). \\[2mm]
In the proof of iii) and in the proof of (\ref{b_induction}) we use the
following identity which can be verified by direct calculation:
\begin{lemma}
\label{rearrange_inverses}
If $1+\alpha\beta$ is invertible then 
$1-\beta(1+\alpha\beta)^{-1}\alpha=(1+\beta\alpha)^{-1}$.\hspace*{\fill}\qed
\end{lemma}

\noindent iii) If $\alpha=\left(\begin{matrix}
1 & u \\
0 & 1    
\end{matrix}\right)$ then setting $\gamma=1+ul'$, we have
\begin{align*}
D(\alpha\beta)&=D\left(
\left(\begin{matrix}
1 & u \\
0 & 1    
\end{matrix}\right)
\left(\begin{matrix}
1   &  0 \\
l'  &  1
\end{matrix}\right)
\left(\begin{matrix}
d'_1 & 0 \\
0        & d'_2
\end{matrix}\right)
\left(\begin{matrix}
1 & u' \\
0 & 1    
\end{matrix}\right)
\right) \\
&=
D\left(\left(\begin{matrix}
1  &  0 \\
l'\gamma^{-1}  &  1
\end{matrix}\right)
\left(\begin{matrix}
\gamma d'_1 & 0 \\
0          & (1-l'\gamma^{-1}u)d'_2
\end{matrix}\right)
\left(\begin{matrix}
1 & (\gamma d'_1)^{-1}ud'_2+u' \\
0 & 1    
\end{matrix}\right)\right). \\
&= \gamma d'_1D((1-l'\gamma^{-1}u)d'_2). \\
&= (1+ul')d'_1D\left((1+l'u)^{-1}d'_2\right) \quad \mbox{(by lemma~\ref{rearrange_inverses})} 
\end{align*}
It follows by $\mfa(n-1)$ and $\mfb(1,n-1)$ that
\begin{equation*}
D(\alpha\beta) = (1+ul')d'_1(1+ul')^{-1}D(d'_2) = d'_1D(d'_2) =
D(\alpha)D(\beta)
\end{equation*}
This completes the proof of (\ref{a_induction}). \\[2mm] \noindent
{\em Proof of (\ref{b_induction}).} Suppose $\alpha:B^n\to B^m$,
$\beta:B^m\to B^n$ and either $\epsilon(\alpha)=0$ or
$\epsilon(\beta)=0$. Write
$\alpha=(\begin{matrix}
\alpha_1 & \alpha_2
\end{matrix})$ where $\alpha_1:B\to B^m$ and
$\alpha_2:B^{n-1}\to B^m$ and write
$\beta=\left(\begin{matrix}
\beta_1 \\
\beta_2
\end{matrix}\right)$ where $\beta_1:B^m\to B$ and $\beta_2:B^m\to B^{n-1}$.
Now 
\begin{align*}
1+\beta\alpha&=
1+
\left(\begin{matrix}\beta_1 \\
\beta_2
\end{matrix}\right)
\bigl(\begin{matrix}
\alpha_1 & \alpha_2
\end{matrix}\bigr)=
\left(\begin{matrix}
1+\beta_1\alpha_1 & \beta_1\alpha_2 \\
\beta_2\alpha_1 & 1+\beta_2\alpha_2
\end{matrix}\right) \\
\intertext{so by definition~\ref{Define_D} we have}
D(1+\beta\alpha)&=(1+\beta_1\alpha_1)D(1+\beta_2\alpha_2-\beta_2\alpha_1(1+\beta_1\alpha_1)^{-1}\beta_1\alpha_2)
\\
&=
(1+\beta_1\alpha_1)D(1+\beta_2(1-\alpha_1(1+\beta_1\alpha_1)^{-1}\beta_1)\alpha_2) \\
&=(1+\beta_1\alpha_1)D(1+\beta_2(1+\alpha_1\beta_1)^{-1}\alpha_2) \quad
\mbox{(by lemma~\ref{rearrange_inverses})} \\
&= D(1+\alpha_1\beta_1)D(1+(1+\alpha_1\beta_1)^{-1}\alpha_2\beta_2) \\
&\mbox{(using $\mfb(1,m)$ and $\mfb(n-1,m)$)} \\
&= D((1+\alpha_1\beta_1)(1+(1+\alpha_1\beta_1)^{-1}\alpha_2\beta_2))
\quad \mbox{(by $\mfa(m)$)} \\
&=D(1+\alpha\beta)
\end{align*}
This completes the proof of proposition~\ref{D_respects_commutation}.
%\hspace*{\fill}\qed
\end{proof}
Thus $D:\tK_1(B)\to \epsilon^{-1}(1)/C_0$ is
well-defined and the proof of theorem~\ref{identify_K1} is also complete. \hspace*{\fill}\qed
%
%\begin{ack}
\vspace*{3ex}

I am grateful to my former PhD adviser Andrew Ranicki who
planted in my mind the questions from which this work germinated. I
also thank Jonathan Rosenberg for bringing to my attention
lemma~\ref{adding_commutation} and Pere Ara and Paul Cohn 
for helpful comments on an earlier draft of the paper.
%\end{ack}
%
%The following command reduces the spacing between bibliographic entries.
\providecommand{\bysame}{\leavevmode\hbox to3em{\hrulefill}\thinspace}
%

%\bibliographystyle{plain}
%\bibliography{../biblio/abbrv}

\begin{thebibliography}{10}

\bibitem{Alm73}
G.~Almkvist.
\newblock Endomorphisms of finitely generated projective modules over a
  commutative ring.
\newblock {\em Arkiv f\"or Matematik}, 11:263--301, 1973.

\bibitem{Alm74}
G.~Almkvist.
\newblock The {Grothendieck} ring of the category of endomorphisms.
\newblock {\em Journal of Algebra}, 28:375--388, 1974.

\bibitem{Alm78}
G.~Almkvist.
\newblock {$K$-theory} of endomorphisms.
\newblock {\em Journal of Algebra}, 55(2):308--340, 1978.

\bibitem{Ara03}
P.~Ara.
\newblock Finitely presented modules over {Leavitt} algebras.
\newblock Preprint.

\bibitem{AGP02}
P.~Ara, K.~R. Goodearl, and E.~Pardo.
\newblock {$K_0$} of purely infinite simple regular rings.
\newblock {\em $K$-Theory}, 26:69--100, 2002.

\bibitem{Bas68}
H.~Bass.
\newblock {\em {Algebraic {$K$-theory}}}.
\newblock W. A. Benjamin Inc, New York-Amsterdam, 1968.

\bibitem{Bot88}
R.~Bott.
\newblock Morse theory indomitable.
\newblock {\em I.H.E.S. Publications Math\'ematiques}, (68):99--114, 1988.

\bibitem{Coh71}
P.~M. Cohn.
\newblock {\em {Free Rings and their Relations}}.
\newblock London Mathematical Society Monographs, 2. Academic Press, London,
  1971.

\bibitem{Coh85}
P.~M. Cohn.
\newblock {\em {Free Rings and their Relations}}.
\newblock London Mathematical Society Monographs, 19. Academic Press, London,
  2nd edition, 1985.

\bibitem{Die43}
J.~Dieudonn\'e.
\newblock Les d\'eterminants sur un corps non commutatif.
\newblock {\em Bulletin de la Soci\'et\'e Math\'ematique de France}, 71:27--45,
  1943.

\bibitem{Dra83}
P.~K. Draxl.
\newblock {\em {Skew Fields}}.
\newblock London Mathematical Society Lecture Note Series, 81. Cambridge
  University Press, 1983.

\bibitem{Gra77}
D.~R. Grayson.
\newblock The {K}-theory of endomorphisms.
\newblock {\em Journal of Algebra}, 48:439--446, 1977.

\bibitem{Har77}
R.~Hartshorne.
\newblock {\em Algebraic Geometry}.
\newblock Number~52 in Graduate Texts in Mathematics. Springer, New
  York-Heidelberg, 1977.

\bibitem{Kli62}
W.~Klingenberg.
\newblock Die struktur der linearen gruppe \"uber einem nichtkommutativen
  lokalen ring.
\newblock {\em Arch. Math.}, 13:73--81, 1962.

\bibitem{Lam01}
T.~Y. Lam.
\newblock {\em A First Course in Noncommutative Rings}.
\newblock Springer, New-York, 2nd edition, 2001.

\bibitem{MenMon84}
P.~Menal and J.~Moncasi.
\newblock {$K_1$} of von {Neumann} regular rings.
\newblock {\em Journal of Pure and Applied Algebra}, 1984.

\bibitem{Mil63}
J.~Milnor.
\newblock {\em Morse theory}.
\newblock Number~51 in Annals of Mathematics Studies. Princeton University
  Press, N.J., 1963.

\bibitem{Mil65}
J.~Milnor.
\newblock {\em Lectures on the $h$-cobordism theorem}.
\newblock Princeton University Press, Princeton, N.J., 1965.
\newblock Notes by L.Siebenmann and J.Sondow.

\bibitem{NeeRan01}
A.~Neeman and A.~A. Ranicki.
\newblock Noncommutative localization and chain complexes {I}. {Algebraic}
  {$K$-} and {$L$-} theory.
\newblock arXiv:math.RA/0109118.

\bibitem{Nov81}
S.P. Novikov.
\newblock Multivalued functions and functionals. {An} analogue of {Morse}
  theory.
\newblock {\em Dokl.Akad.Nauk.SSSR}, 260(1):31--35, 1981.
\newblock English translation: Soviet Math.Dokl., 24 (1981) no.2 222-226.

\bibitem{Nov82}
S.P. Novikov.
\newblock The {Hamiltonian} formalism and a multivalued analogue of {Morse}
  theory.
\newblock {\em Uspeki Mat.}, 37(5):3--49, 1982.
\newblock English translation: Russian Math.Surveys 37 (1982) 1-56.

\bibitem{PajRan00}
A.~Pajitnov and A.~A. Ranicki.
\newblock {The Whitehead group of the Novikov ring}.
\newblock {\em $K$-theory}, 21(4):325--365, 2000.
\newblock arXiv:math.AT/0012031.

\bibitem{Paj95}
A.~V. Pajitnov.
\newblock On the {Novikov} complex for rational {Morse} forms.
\newblock {\em Annales de la Facult\'e des Sciences de Toulouse},
  IV(2):297--338, 1995.

\bibitem{Paj99}
A.~V. Pajitnov.
\newblock Simple homotopy type of the {Novikov} complex, and the {Lefschetz}
  $\zeta$-function of the gradient flow.
\newblock {\em Uspekhi Mat. Nauk}, 54(1(325)):117--170, 1999.
\newblock Translation in Russian Math. Surveys 54(1) 119-169, 1999.
  arXiv:math.dg-ga/9706014.

\bibitem{Paj00}
A.~V. Pajitnov.
\newblock Closed orbits of gradient flows and logarithms of non-abelian {Witt}
  vectors.
\newblock {\em $K$-Theory}, 21(4):301--324, 2000.
\newblock arXiv:math.DG/99808010.

\bibitem{Ran98}
A.~A. Ranicki.
\newblock {\em High-dimensional Knot Theory}.
\newblock Springer, Berlin, 1998.

\bibitem{Ran01}
A.~A. Ranicki.
\newblock Circle valued {Morse} theory and {Novikov} homology.
\newblock In F.~T. Farrell, L.~G\"ottsche, and W.~L\"uck, editors, {\em
  Topology of High-Dimensional Manifolds}, number~2 in ICTP Lecture Notes. The
  Abdus Salam ICTP, Trieste, 2002.
\newblock Proceedings of the Summer School on High-dimensional Manifold
  Topology, 21 May - 8 June 2001. arXiv:math.AT/0111317.

\bibitem{Rev83}
G.~R\'ev\'esz.
\newblock On the abelianized multiplicative group of a universal field of
  fractions.
\newblock {\em Journal of Pure and Applied Algebra}, 27(3):277--297, 1983.

\bibitem{Ros94}
J.~Rosenberg.
\newblock {\em {Algebraic {$K$}-theory and its applications}}.
\newblock Graduate Texts in Mathematics, 147. Springer, New York, 1994.

\bibitem{Scho85}
A.~H. Schofield.
\newblock {\em {Representations of rings over skew fields}}, volume~92 of {\em
  London Mathematical Society Lecture Note Series}.
\newblock Cambridge University Press, 1985.

\bibitem{She01}
D.~Sheiham.
\newblock Non-commutative characteristic polynomials and {Cohn} localization.
\newblock {\em Journal of the London Mathematical Society (2)}, 64(1):13--28,
  2001.
\newblock arXiv:math.RA/0104158.

\bibitem{Sil81}
J.~R. Silvester.
\newblock {\em Introduction to algebraic $K$-theory}.
\newblock Chapman and Hall, London-New York, 1981.

\bibitem{Sri96}
V.~Srinivas.
\newblock {\em Algebraic $K$-theory}, volume~90 of {\em Progress in
  Mathematics}.
\newblock Birkh\"auser, Boston, 1996.

\end{thebibliography}
\vspace*{3ex}
Department of Mathematics \\ University of California,
Riverside \\ CA 92521, USA \\ des@sheiham.com
\end{document}